\documentclass[11pt]{article}%[10pt]{article}%\documentclass
\usepackage{amsmath,amssymb,latexsym,epsf,mathabx}%

\begin{document}

%% ----------------------------------------------------------------------
\newcommand{\nc}{\newcommand}
\def\PP#1#2#3{{\mathrm{Pres}}^{#1}_{#2}{#3}\setcounter{equation}{0}}
\def\ns{$n$-star}\setcounter{equation}{0}
\def\nt{$n$-tilting}\setcounter{equation}{0}
\def\Ht#1#2#3{{{\mathrm{Hom}}_{#1}({#2},{#3})}\setcounter{equation}{0}}
\def\qp#1{{${(#1)}$-quasi-projective}\setcounter{equation}{0}}
\def\mr#1{{{\mathrm{#1}}}\setcounter{equation}{0}}
\def\mc#1{{{\mathcal{#1}}}\setcounter{equation}{0}}
\def\HD{\mr{Hom}_{\mc{D}}}
\def\AdT{\mr{add}_{\mc{D}}}
\def\Kb{\mc{K}^b(\mr{Proj}R)}
\def\kb{\mc{K}^b(\mc{P}_R)}

%%%%%%%%FROM latexexam.tex
%\theoremstyle{definition}
\newtheorem{Th}{Theorem}[section]%[section]%
\newtheorem{Def}[Th]{Definition}
\newtheorem{Lem}[Th]{Lemma}
\newtheorem{Pro}[Th]{Proposition}
\newtheorem{Cor}[Th]{Corollary}
\newtheorem{Rem}[Th]{Remark}
\newtheorem{Exm}[Th]{Example}
\newtheorem{Sc}[Th]{}
\def\Pf#1{{\noindent\bf Proof}.\setcounter{equation}{0}}
\def\>#1{{ $\Rightarrow$ }\setcounter{equation}{0}}
\def\<>#1{{ $\Leftrightarrow$ }\setcounter{equation}{0}}
\def\bskip#1{{ \vskip 20pt }\setcounter{equation}{0}}
\def\sskip#1{{ \vskip 5pt }\setcounter{equation}{0}}
\def\mskip#1{{ \vskip 10pt }\setcounter{equation}{0}}
\def\bg#1{\begin{#1}\setcounter{equation}{0}}
\def\ed#1{\end{#1}\setcounter{equation}{0}}
\def\KET{T^{^F\bot}\setcounter{equation}{0}}
\def\KEC{C^{\bot}\setcounter{equation}{0}}
\def\Db{\mc{D}^b(\mr{mod}R)}

%%%%%%%%%%
%
%\def\jze{{ \begin{pmatrix} 0 & 0 \\ 1 & 0 \end{pmatrix}}\setcounter{equation}{0}}
%\def\hjz#1#2{{ \begin{pmatrix} {#1} & {#2} \end{pmatrix}}\setcounter{equation}{0}}
%\def\ljz#1#2{{  \begin{pmatrix} {#1} \\ {#2} \end{pmatrix}}\setcounter{equation}{0}}
%\def\jz#1#2#3#4{{  \begin{pmatrix} {#1} & {#2} \\ {#3} & {#4} \end{pmatrix}}\setcounter{equation}{0}}

%%%%%%%%%%%%%%%%%%%%%%%%%%%%%%%%%%%%%%%%%%%%%%%%%%%%%%%%%%%%%%%%%%%%%%%%%%%%%%%%%
%**************************±êÌâ¡¢ÕªÒª¡¢·ÖÀàºÅ¡¢¹Ø¼ü×Ö**************************
%
%\title{\bf Relative singularity categories and $\omega$-Gorenstein objects in triangulated categories {\thanks {Supported by the National Science Foundation of
%
\title{\bf Relative singularity categories, Gorenstein objects  and silting theory {\thanks {Supported by the National Science Foundation of China
(Grant Nos. 11171149) and the National Science Foundation for Distinguished Young Scholars of Jiangsu Province (Grant No.BK2012044) and a project funded by the Priority Academic Program Development of Jiangsu Higher Education Institutions}}}
\smallskip
\author{\small {Jiaqun WEI}\\
\small Institute of Mathematics, School of Mathematics Sciences,\\
\small
Nanjing Normal University \\
\small Nanjing 210023, P.R.China\\ \small Email:
weijiaqun@njnu.edu.cn}
\date{}
\maketitle
\baselineskip 16pt%18pt%14pt%15.5pt%\baselineskip  25.5pt %
%%%%%%%\hskip 18pt
%
% Abstract ------------------------------------------------------
%
\begin{abstract}
\vskip 10pt%

We study singularity categories through Gorenstein objects in triangulated categories and silting theory. Let ${\omega}$ be a semi-selforthogonal (or presilting) subcategory of a triangulated category $\mc{T}$. We introduce the notion of $\omega$-Gorenstein objects, which is far extended version of Gorenstein projective modules and Gorenstein injective modules in triangulated categories. We prove that the stable category $\underline{\mc{G}_{\omega}}$, where $\mc{G}_{\omega}$ is the subcategory of all ${\omega}$-Gorenstein objects, is a triangulated category and it is, under some conditions, triangle equivalent to the relative singularity category of $\mc{T}$ with respect to $\omega$. As applications, we obtain the following characterizations of singularity categories which partially extend classic results (usually in the context of Gorenstein rings) to some more general settings. (1) For a ring $R$ of finite Gorenstein global dimension, there are triangle equivalences between $\underline{\mc{G}P}$ (the stable category of Gorenstein projective modules), $\underline{\mc{G}I}$ (the stable category of Gorenstein injective modules), $\underline{\mc{G}_{\mr{Add}M}}$ (where $M$ is any big silting complex in $\mc{D}^b(\mr{Mod}R)$) and the big singularity category $\mc{D}_{\mr{Sg}}(R)$. (2) For a left coherent ring $R$ of finite Gorenstein global dimension, there are triangle equivalences between $\underline{\mc{G}p}$ (the stable category of finitely generated Gorenstein projective modules), $\underline{\mc{G}_M}$ (where $M$ is any silting complex in $\mc{D}^b(R)$) and the singularity category $\mc{D}_{sg}(R)$.

\mskip\

\noindent MSC 2010: %Mathematics Subject Classification:
  Primary 18E30 16E05 Secondary 18G35 16G10

\sskip\

\noindent {\it Keywords}: Gorenstein object, silting complex, singualrity category, triangulated category

\end{abstract}
%\smallskip
%
\vskip 30pt
% ----------------------------------------------------------------------
%% ----------------------------------------------------------------------
%\def\baselinestretch{1}

\section{Introduction}
%{\noindent \Large\bf Introduction}

%
%\bskip\
%
%
\hskip 18pt
%%% ----------------------------------------------------------------------
%\vskip 30pt
%---------------------------------------------------------------------------------------------
%\bskip\

\mskip\

Singularity categories (with another name stable derived categories [K]) are important in the study of  varieties and rings of infinite global dimension.
For an algebraic variety $X$, the singularity category is given by a Verdier's quotient triangulated category $\mc{D}_{sg}(X):=\mc{D}^b(\mr{coh}X)/\mr{perf}X$, where $\mc{D}^b(\mr{coh}X)$ is bounded derived category of coherent sheaves on $X$ and  $\mr{perf}X$ is its full subcategory of perfect complexes [O]. Note that $\mc{D}_{sg}(X)=0$ if and only if $X$ is smooth. In the study of representation theory of algebras, singularity categories are defined as the Verdier's quotient triangulated categories $\mc{D}_{sg}(R):=\mc{D}^b(R)/\mc{K}^b(\mc{P})$, where $\mc{D}^b(R)$ is the bounded derived category of finitely presented modules over a (left) coherent ring $R$ and  $\mc{K}^b(\mc{P})$ is the bounded homotopy category of finitely generated projective $R$-modules [Bu, C]. Note that the singularity category $\mc{D}_{sg}(R)$ is 0 if and only if $R$ has finite global dimension.

Similar quotient triangulated categories were also studied by several people, such as $\mc{D}_{\mr{Sg}}(R):=\mc{D}^b(\mr{Mod}R)/\mc{K}^b(\mr{Proj})$ [Be1], $\mc{D}^b(\mr{Mod}R)/\mc{K}^b(\mr{Inj})$ [Be1],  $\mc{D}^b(R)/\mc{K}^b(\mc{I})$ [H2], $\mc{D}^b(R)/\mc{K}^b(\mr{add}T)$ [CZ], where $\mr{Mod}R$ ($\mr{Proj}$, $\mr{Inj}$, $\mc{I}$, resp.) denotes over a ring $R$ the category of all (projective, injective, finitely generated injective, resp.) $R$-modules. In [C], a more general notion of relative singularity categories $\mc{D}_{\omega}(\mc{A}):=\mc{D}^b(\mc{A})/\mc{K}^b(\mc{\omega})$ was prompted, where $\mc{A}$ is an abelian category and $\omega$ is a selforthogonal additive subcategory of $\mc{A}$. The notion unified the previous variant singularity categories.

Note that Gorenstein projective modules (or Cohen-Macauley objects) and Gorenstein injective modules over Gorenstein rings play important roles in the study of the above relative singularity categories [Be1, Bu, H2, C], etc. The theory of such Gorenstein modules  stems from  Auslander's notion of modules of G-dimension zero [AB] and later were extensively studied and developed by Avramov-Foxby [AF], Buchwitz [Bu], Enochs-Jenda [EJ], etc. and their schools in variant ways. It is known that Gorenstein global dimension of a ring is well-defined, and it coincides with the supersum of all Gorenstein projective dimensions of modules as well as the supersum of all Gorenstein injective dimensions of modules [Be1, BM].
An important class of rings with finite Gorenstein global dimension are Gorenstein rings, that is, left and right noetherian rings with finite left and right self-injective dimensions. However, there are also non-noetherian rings with finite Gorenstein global dimension [MT].

In this paper, we aim to study singularity categories through silting theory. Recall that silting theory was first introduced by Keller and Vossieck [KV] and then recognized by Aihara and Iyama [AI] from the point of mutations and by the auther [W2] (with the name semi-tilting complexes) from the point of compare between tilting complexes and tilting modules. In particular, from the point of [W2], it is reasonable to say that the role that silting complexes play in derived categories is similar as that tilting modules play in module categories, while the role that tilting complexes play in derived categories is that progenerators play in module categories.

%
%extend Gorenstein objects and singularity categories  to general triangulated categories and study their relations, through silting theory.

The paper is organized as follows. After the introduction, we present in Section 2 a new extension of Gorenstein objects in triangulated categories from the point of tilting theory (compare [AS]) and introduce the relative singularity category $\mc{T}_{\omega}$ for a triangulated category $\mc{T}$ with a semi-selforthogonal subcategory $\omega$, by employing methods and techniques in the study of silting theory in triangulated categories [W2]. This generalizes the above-mentioned variant singularity categories.

The following is our main result which characterizes the relative singularity category $\mc{T}_{\omega}$ as the stable category of $\omega$-Gorenstein objects in $\mc{T}$ under some conditions, see Theorem \ref{T} for details.

%\baselineskip 15pt
%\vskip 10pt
\bg{quote}%
\noindent {\bf Theorem A} Let $\mc{T}$ be a triangulated category and $\omega$ be a semi-selforthogonal subcategory of $\mc{T}$. If $\mc{T}=\langle {_{\omega}\mc{X}}\rangle=\langle {\mc{X}_{\omega}}\rangle$, then there is a triangle equivalence between $\underline{\mc{G}_{\omega}}$ and the relative singularity $\mc{T}_{\omega}$.
\ed{quote}%

 Applications of Theorem A are given in Section 3, where we obtain the following results which extend classical results [Be1, Bu, C], etc. to more general settings (Theorem B (2)) and provide characterizations of singularity categories in terms of silting complexes, see Section 3 for details. Readers who are interesting in algebraic geometry could deduce some other corollaries.

%\vskip 10pt
%\baselineskip 15pt
\bg{quotation}%
\noindent {\bf Theorem B} (1) Let $R$ be a ring with finite Gorenstein global dimension. Then there are triangle equivalences between the following triangulated categories.
%\bg{quotation}%

 (i) $\underline{\mc{G}P}$ (the stable category of Gorenstein projective modules),

 (ii) $\underline{\mc{G}I}$ (the stable category of Gorenstein injective modules),

 (iii) $\underline{\mc{G}_{\mr{Add}M}}$ (the stable category of ($\mr{Add}M$)-Gorenstein objects in $\mc{D}^b(\mr{Mod}R)$, where $M$ is any big silting complex in $\mc{D}^b(\mr{Mod}R)$) and,

 (iv) $\mc{D}_{\mr{Sg}}(R)$ (the big singularity category $\mc{D}^b(\mr{Mod}R)/\mc{K}^b(\mr{Proj})$).
%\ed{quotation}%

\noindent (2) Let $R$ be a left coherent ring with finite Gorenstein global dimension. Then there are triangle equivalences between the following triangulated categories.
%\bg{quotation}%

 (i) $\underline{\mc{G}p}$ (the stable category of finitely generated Gorenstein projective modules),

 (ii) $\underline{\mc{G}_M}$ (the stable category of ($\mr{add}M$)-Gorenstein objects in $\mc{D}^b(R)$, where $M$ is any silting complex in $\mc{D}^b(R)$) and,

 (iii) $\mc{D}_{sg}(R)$ (the singularity category $\mc{D}^b(R)/\mc{K}^b(\mc{P})$).
%\ed{quotation}%

\ed{quotation}%

 Complexes in the paper are always cochain ones and subcategories
are always full.
 We
denote by $fg: L\to N$ the composition of two homomorphism $f:L\to
M$ and $g:M\to N$. For simple, we use $\mr{Hom}_{\mc{D}}(-,-)$ instead of $\mr{Hom}_{\mc{D}^b(\mr{Mod}R)}(-,-)$.

For basic knowledge on triangulated categories,  derived
categories and the tilting theory, we refer to [{H1}].% and the
%Handbook of tilting theory [{AHKb}].

%
\bskip\
\vskip 30pt
% ----------------------------------------------------------------------
%% ----------------------------------------------------------------------
%\def\baselinestretch{1}

\section{Main results}
%%
%%
%%%% ----------------------------------------------------------------------
%%
%\mskip\
%

%%%%%%%%%%%%%%
\def\HD{\mr{Hom}_{\mc{D}}}
\def\HT{\mr{Hom}_{\mc{T}}}
%%%%%%%%%%%%%%

%%%%%%%%%%%%%%%%%%%%%%%%%%%%%%%%%%%%%%%%%%%%%%%%%%%%%%%%%%%%%%%%%%%%%%%%%%%%%%%%%%
\hskip 18pt
%
%Hom-finite

Let $\mc{T}$ be a triangulated category with the associated shift functor [1].

Let $\mc{C}$ be a subcategory in $\mc{T}$. Similarly as the notions in [AR], we denote $\widehat{\mc{C}}$ to be the class of all objects $T$ satisfying that there are triangles $T_{i+1}\to C_i\to T_i\to$ for some $r\ge 0$ and all $0\le i\le r$, such that $T_0=T$ and $T_{r+1}=0$, $C_i\in \mc{C}$ for each $i$. Dually, we denote $\widecheck{\mc{C}}$ to be the class of all objects $T$ satisfying that there are triangles $T_{i}\to C_i\to T_{i+1}\to$ for some $r\ge 0$ and all $0\le i\le r$, such that $T_0=T$ and $T_{r+1}=0$, $C_i\in \mc{C}$ for each $i$. As usual, we denote by $\langle\mc{C}\rangle$ to be the smallest triangulated subcategory containing $\mc{C}$.

Let $\mc{C}$ be a subcategory  containing  $0$. The  subcategory $\mc{C}$ is
extension-closed if for any triangle $U\to V\to W\to$ with
$U,W\in\mc{C}$, it holds that $V\in \mc{C}$.  It is resolving
(resp., coresolving) if it is further closed under the functor
$[-1]$ (resp., $[1]$). Note that $\mc{C}$ is resolving (resp.,
coresolving) if and only if, for any triangle $U\to V\to W\to$
(resp.,  $W\to V\to U\to$) in $\mc{T}$ with $W\in\mc{C}$, one has
that ``\ $U\in\mc{C}\Leftrightarrow V\in\mc{C}$\ ''.

%%%%%%%%%%%%------------------------------------------------------------------------------------------------------------------------------
\vskip 10pt
The following result is useful.

\bg{Lem}\label{Rl}% resolving lemma
Let $\mc{C}$ be a subcategory of $\mc{T}$ containing $0$. Denote that $\mc{C}[-]=\{X\in\mc{T}\ |\ X\simeq C[n]$ for some $C\in\mc{C}$ and some $n\le 0\}$. The notion $\mc{C}[+]$ is defined dually.

$(1)$ If $\mc{C}$ is resolving, then $\langle \mc{C}\rangle=\mc{C}[+]=\widehat{\mc{C}}$.

$(2)$ If $\mc{C}$ is coresolving, then $\langle \mc{C}\rangle=\mc{C}[-]=\widecheck{\mc{C}}$.
\ed{Lem}

\Pf. (1) We firstly prove that $\langle C\rangle=\mc{C}[+]$ and then prove that they coincide with $\widehat{\mc{C}}$.

 Clearly $\mc{C}\subseteq\mc{C}[+]\subseteq\langle C\rangle$ by the definitions. It is easy to see that $\mc{C}[+]$ is closed under $[1], [-1]$. Now, take any triangle $X \to Y\to Z\to$ with $X,Z\in\mc{C}[+]$. Note that there is some integer $n$ such that $X,Z\in\mc{C}[n]$, i.e., $X[-n],Z[-n]\in\mc{C}$. If $n\le 0$, we already have that $Y\in\mc{C}\subseteq \mc{C}[+]$ since $\mc{C}$ is resolving. If $n>0$, then we have the triangle $X[-n]\to Y[-n]\to Z[-n]\to$. But $X[-n],Z[-n]\in\mc{C}$ implies that $Y[-n]\in\mc{C}$ since $\mc{C}$ is closed under extensions. It follows that $Y=(Y[-n])[n]\in\mc{C}[+]$. Altogether we obtain that $\mc{C}[+]$ is a triangulated subcategory of $\mc{T}$. Then we see that $\langle C\rangle=\mc{C}[+]$.

 It is obvious that $\widehat{\mc{C}}\subseteq\langle C\rangle$. On the other hand, take any $X=C[n]$ for some $C\in\mc{C}$. If $n\le 0$, we already have that $X\in\mc{C}\subseteq \widehat{\mc{C}}$, since $\mc{C}$ is resolving. Otherwise, $n>0$ and since $0\in\mc{C}$, we easily obtain $X=C[n]\in \widecheck{\mc{C}}$. Then we conclude that $\langle C\rangle=\mc{C}[+]=\widecheck{\mc{C}}$ by the definition.

 (2) The proof is dual to (1).
\hfill$\Box$

%%%%%%%%%%%%------------------------------------------------------------------------------------------------------------------------------
\vskip 10pt

A subcategory ${\omega}\subseteq \mc{T}$ is {\it{semi-selforthogonal
}} ({\it{selforthogonal}}
 resp.) if $\HT{(M, M'[i])} =0$ for all $M,M'\in\omega$ and $i>0$ ($i\neq 0$, resp.). As examples, the subcategory of projective modules over a ring $R$ is selforthogonal in the derived category $\mc{D}^b(\mr{Mod}R)$, as well as the subcategory of injective modules. Clearly, $\omega[i]$ is also semi-selforthogonal for any integer $i$, in case $\omega$ is semi-selforthogonal. To distinguish, we call an object $M\in\mc{T}$  {\it{semi-selforthogonal}} ({\it{product-semi-selforthogonal}}, {\it{coproduct-semi-selforthogonal}}, resp.) if the subcategory $\{M\}$ ($\{M^X\ |\ X$ is an index set$\}$, $\{M^{(X)}\ |\ X$ is an index set$\}$, resp.) is semi-selforthogonal. Similarly, one also have notions of {\it{selforthogonal}}, {\it{product-selforthogonal}}, {\it{coproduct-selforthogonal}} objects. We note that semi-selforthogonal subcategories (objects, resp.) are also called presilting subcategories (objects, resp.) in some papers.

From now on, we fix ${\omega}$ to be semi-selforthogonal. If $\omega=\{M\}$, we often replace $\omega$ with $M$ in the correspondent notions.

As usual, we denote by ${^{\perp}{\omega}}$ the class of all objects $N$ such that $\HT (N,{\omega}[i])=0$ for all $i>0$ and by  ${\omega}^{\perp}$ the class of all objects $N$ such that $\HT ({\omega},N[i])=0$ for all $i>0$. Also we usually denote by $\mr{add}{\omega}$ the class of all direct summands of finite direct sums of copies of objects in ${\omega}$. Obviously, $\mr{add}{\omega}\subseteq {^{\perp}{\omega}}\bigcap {\omega}^{\perp}$ holds in our setting.

% {\it\textbf{Auslander class}}
We define the {\it{Auslander class}} associated to ${\omega}$, denoted by $\mc{X}_{\omega}$, to be the class of all objects $T$ satisfying that there are triangles $T_i\to M_i\to T_{i+1}\to$ for all $i\ge 0$, such that $T_0=T$, $M_i\in \mr{add}{\omega}$ and $T_i\in {^{\perp}{\omega}}$, for each $i\ge 0$. The notion of the Auslander class follows from [{AR}, Section 5], where it was defined for a selforthogonal subcategory of modules $\omega$. We have a dual notion called {\it{ co-Auslander class}} associated to ${\omega}$, denoted by ${_{\omega}\mc{X}}$, consisting of all objects $T$ satisfying that there are triangles $T_{i+1}\to M_i\to T_i\to$ for all $i\ge 0$, such that $T_0=T$, $M_i\in \mr{add}{\omega}$ and $T_i\in {{\omega}^{\perp}}$, for each $i\ge 0$.

The following result collects some basic properties on subcategories associated with the semi-selforthogonal subcategory ${\omega}$.

%%%%%%%%%%%%%%%%%%%%%%%%%%%%%%%%%%%%%%%%%%%%%%%%%%%%%%%%%%%%%%%%%%%%%%%%%%%%%%%%%%
%\hskip 18pt
%
\bg{Lem}\label{Bl}%Basic lemma for M
Let ${\omega}$ be semi-selforthogonal. Then

$(1)$  $\widecheck{\mr{add}{\omega}} \subseteq \mc{X}_{\omega}\subseteq {^{\perp}{\omega}}$ and all three classes are closed under taking extensions, $[-1]$, direct summands and finite direct sums, i.e., they are resolving.

$(2)$  $\widehat{\mr{add}{\omega}} \subseteq {_{\omega}\mc{X}}\subseteq {{\omega}^{\perp}}$ and all three classes are closed under taking extensions, $[1]$, direct summands and finite direct sums, i.e., they are coresolving.

$(3)$ For any $X\in\widecheck{_{\omega}\mc{X}}$, there are triangles $U\to V\to X\to$ and $X\to U'\to V'\to$ with $U,U'\in{_{\omega}\mc{X}}$ and $V,V'\in\widecheck{\mr{add}{\omega}}$.

$(4)$ For any $X\in\widehat{\mc{X}_{\omega}}$, there are triangles $X\to V\to U\to$ and $V'\to U'\to X\to$ with $U,U'\in\mc{X}_{\omega}$ and $V,V'\in\widehat{\mr{add}{\omega}}$.

\ed{Lem}

\Pf. (1) and (2). These two results were already proved in [W2, Section 2] in case $\omega={\mr{Add}M}$ for a coproduct-semi-selforthogonal object $M$. Their proofs can be transferred completely.

(3) By imitating the proofs of [W2, Lemma 2.7 and Corollary 2.8], we obtain that the first triangle exists, i.e., for any $X\in\widecheck{_{\omega}\mc{X}}$, there is a triangle $U\to V\to X\to$ with $U\in{_{\omega}\mc{X}}$ and $V\in\widecheck{\mr{add}{\omega}}$. Note that there is a triangle $V\to M_V\to V'\to$ with $M_V\in\mr{add}{\omega}$ and $V'\in\widecheck{\mr{add}{\omega}}$, as $V\in\widecheck{\mr{add}{\omega}}$. Thus, we have the following commutative diagram of triangles.

\mskip\

 \setlength{\unitlength}{0.09in}
 \begin{picture}(50,18)

%                 \put(18,3.4){\vector(0,-1){2}}
                 \put(27,3.4){\vector(0,-1){2}}
                 \put(35,3.4){\vector(0,-1){2}}

% \put(18,5){\makebox(0,0)[c]{$U$}}
%                             \put(21,5){\vector(1,0){2}}
 \put(27,5){\makebox(0,0)[c]{$V'$}}
                             \put(30,4.9){\line(1,0){2}}
                             \put(30,5.2){\line(1,0){2}}
 \put(35,5){\makebox(0,0)[c]{$V'$}}

%                 \put(18,9){\vector(0,-1){2}}
                 \put(27,9){\vector(0,-1){2}}
                 \put(35,9){\vector(0,-1){2}}
%                 \put(34.5,9){\line(0,-1){2}}
%                 \put(35,9){\line(0,-1){2}}

 \put(18,11){\makebox(0,0)[c]{$U$}}
                             \put(21,11){\vector(1,0){2}}
 \put(27,11){\makebox(0,0)[c]{$M_V$}}
                             \put(30,11){\vector(1,0){2}}
 \put(35,11){\makebox(0,0)[c]{$U'$}}
                             \put(37,11){\vector(1,0){2}}

                 \put(17.9,14.5){\line(0,-1){2}}
                 \put(18.2,14.5){\line(0,-1){2}}
                 \put(27,14.5){\vector(0,-1){2}}
                 \put(35,14.5){\vector(0,-1){2}}

 \put(18,16){\makebox(0,0)[c]{$U$}}
                              \put(21,16){\vector(1,0){2}}
%                             \put(21,16){\line(1,0){2}}
%                             \put(21,16.5){\line(1,0){2}}
 \put(27,16){\makebox(0,0)[c]{$V$}}
                              \put(30,16){\vector(1,0){2}}
 \put(35,16){\makebox(0,0)[c]{$X$}}
                              \put(37,16){\vector(1,0){2}}

\end{picture}
\sskip\

Since $U,M_V\in {_{\omega}\mc{X}}$ and ${_{\omega}\mc{X}}$ is coresolving by (2), we see that $U'\in {_{\omega}\mc{X}}$ too. It follows that the right column is just the second desired triangle in (3).

(4) The proof can be obtained by the dual statement of the proof of (3).
\hfill $\Box$

%%%%%%%%%%%%%%%%%%%%%%%%%%%%%%%%%%%%%%%%%%%%%%%%%%%%%%%%%%%%%%%%%%%%%%%%%%%%%%%%%%
\vskip 10pt

For later use, we include the following lemma.

\bg{Lem}\label{Tl}% Triangle lemma
Let ${\omega}\subseteq\mc{T}$ be semi-selforthogonal. Then

$(1)$ For a triangle $T\to T'\to T''\to$ in $\widecheck{_{\omega}\mc{X}}$, there is the following commutative diagram of triangles, where  $U_T,U_{T'},U_{T''}\in{_{\omega}\mc{X}}$ and $V_T,V_{T'},V_{T''}\in\widecheck{\mr{add}{\omega}}$.

\mskip\

 \setlength{\unitlength}{0.09in}
 \begin{picture}(50,18)

                 \put(18,3.4){\vector(0,-1){2}}
                 \put(27,3.4){\vector(0,-1){2}}
                 \put(35,3.4){\vector(0,-1){2}}

 \put(18,5){\makebox(0,0)[c]{$T$}}
                             \put(21,5){\vector(1,0){2}}
 \put(27,5){\makebox(0,0)[c]{$T'$}}
                             \put(30,5){\vector(1,0){2}}
 \put(35,5){\makebox(0,0)[c]{$T''$}}
                             \put(37,5){\vector(1,0){2}}

                 \put(18,9){\vector(0,-1){2}}
                 \put(27,9){\vector(0,-1){2}}
                 \put(35,9){\vector(0,-1){2}}
%                 \put(34.5,9){\line(0,-1){2}}
%                 \put(35,9){\line(0,-1){2}}

 \put(18,11){\makebox(0,0)[c]{$V_T$}}
                             \put(21,11){\vector(1,0){2}}
 \put(27,11){\makebox(0,0)[c]{$V_{T'}$}}
                             \put(30,11){\vector(1,0){2}}
 \put(35,11){\makebox(0,0)[c]{$V_{T''}$}}
                             \put(37,11){\vector(1,0){2}}

                 \put(18,14.5){\vector(0,-1){2}}
                 \put(27,14.5){\vector(0,-1){2}}
                 \put(35,14.5){\vector(0,-1){2}}

 \put(18,16){\makebox(0,0)[c]{$U_T$}}
                              \put(21,16){\vector(1,0){2}}
%                             \put(21,16){\line(1,0){2}}
%                             \put(21,16.5){\line(1,0){2}}
 \put(27,16){\makebox(0,0)[c]{$U_{T'}$}}
                              \put(30,16){\vector(1,0){2}}
 \put(35,16){\makebox(0,0)[c]{$U_{T''}$}}
                              \put(37,16){\vector(1,0){2}}

\end{picture}

$(2)$  For a triangle $T\to T'\to T''\to$ in $\widehat{_{\omega}\mc{X}}$, there is the following commutative diagram of triangles, where $Y_T,Y_{T'},Y_{T''}\in{\mc{X}_{\omega}}$ and $X_T,X_{T'},X_{T''}\in\widehat{\mr{add}{\omega}}$.

 \setlength{\unitlength}{0.09in}
 \begin{picture}(50,20)

                 \put(18,3.4){\vector(0,-1){2}}
                 \put(27,3.4){\vector(0,-1){2}}
                 \put(35,3.4){\vector(0,-1){2}}

 \put(18,5){\makebox(0,0)[c]{$T$}}
                             \put(21,5){\vector(1,0){2}}
 \put(27,5){\makebox(0,0)[c]{${T'}$}}
                             \put(30,5){\vector(1,0){2}}
 \put(35,5){\makebox(0,0)[c]{${T''}$}}
                             \put(37,5){\vector(1,0){2}}

                 \put(18,9){\vector(0,-1){2}}
                 \put(27,9){\vector(0,-1){2}}
                 \put(35,9){\vector(0,-1){2}}
%                 \put(34.5,9){\line(0,-1){2}}
%                 \put(35,9){\line(0,-1){2}}

 \put(18,11){\makebox(0,0)[c]{$Y_T$}}
                             \put(21,11){\vector(1,0){2}}
 \put(27,11){\makebox(0,0)[c]{$Y_{T'}$}}
                             \put(30,11){\vector(1,0){2}}
 \put(35,11){\makebox(0,0)[c]{$Y_{T''}$}}
                             \put(37,11){\vector(1,0){2}}

                 \put(18,14.5){\vector(0,-1){2}}
                 \put(27,14.5){\vector(0,-1){2}}
                 \put(35,14.5){\vector(0,-1){2}}

 \put(18,16){\makebox(0,0)[c]{$X_{T}$}}
                              \put(21,16){\vector(1,0){2}}
%                             \put(21,16){\line(1,0){2}}
%                             \put(21,16.5){\line(1,0){2}}
 \put(27,16){\makebox(0,0)[c]{$X_{T'}$}}
                              \put(30,16){\vector(1,0){2}}
 \put(35,16){\makebox(0,0)[c]{$X_{T''}$}}
                              \put(37,16){\vector(1,0){2}}

\end{picture}
\ed{Lem}

\Pf. (1) Let $U_T\to V_T\to T\to$ and $U_{T''}\to V_{T''}\to T''\to$ with $U_T,U_{T''}\in{_{\omega}\mc{X}}$ and $V_T,V_{T''}\in\widecheck{\mr{add}{\omega}}$ be triangles given by Lemma \ref{Bl}(3). Consider the following diagram.

 \setlength{\unitlength}{0.09in}
 \begin{picture}(50,20)

                 \put(9,3.4){\vector(0,-1){2}}
                 \put(18,3.4){\vector(0,-1){2}}
                 \put(27,3.4){\vector(0,-1){2}}
                 \put(35,3.4){\vector(0,-1){2}}

 \put(9,5){\makebox(0,0)[c]{$T''[-1]$}}
                             \put(13,5){\vector(1,0){2}}
 \put(18,5){\makebox(0,0)[c]{$T$}}
                             \put(21,5){\vector(1,0){2}}
 \put(27,5){\makebox(0,0)[c]{$T'$}}
                             \put(30,5){\vector(1,0){2}}
 \put(35,5){\makebox(0,0)[c]{$T''$}}
                             \put(37,5){\vector(1,0){2}}

                 \put(9,9){\vector(0,-1){2}}
                 \put(18,9){\vector(0,-1){2}}
                 \put(27,9){\vector(0,-1){2}}
                 \put(35,9){\vector(0,-1){2}}
%                 \put(34.5,9){\line(0,-1){2}}
%                 \put(35,9){\line(0,-1){2}}

 \put(9,11){\makebox(0,0)[c]{$V_{T''}[-1]$}}
      \put(14,11.5){\makebox(0,0)[c]{$\stackrel{f}{\dashrightarrow}$}}
%                             \put(13,11){\vector(1,0){2}}
%                             \put(14,12){\makebox(0,0)[c]{$f$}}
 \put(18,11){\makebox(0,0)[c]{$V_T$}}
      \put(22,11){\makebox(0,0)[c]{${\dashrightarrow}$}}
%                             \put(21,11){\vector(1,0){2}}
 \put(27,11){\makebox(0,0)[c]{$V_{T'}$}}
      \put(31,11){\makebox(0,0)[c]{${\dashrightarrow}$}}
%                             \put(30,11){\vector(1,0){2}}
 \put(35,11){\makebox(0,0)[c]{$V_{T''}$}}
                             \put(37,11){\vector(1,0){2}}

                 \put(9,14.5){\vector(0,-1){2}}
                 \put(18,14.5){\vector(0,-1){2}}
                 \put(27,14.5){\vector(0,-1){2}}
                 \put(35,14.5){\vector(0,-1){2}}

 \put(9,16){\makebox(0,0)[c]{$U_{T''}[-1]$}}
      \put(14,16){\makebox(0,0)[c]{${\dashrightarrow}$}}
%                              \put(13,16){\vector(1,0){2}}
 \put(18,16){\makebox(0,0)[c]{$U_T$}}
      \put(22,16){\makebox(0,0)[c]{${\dashrightarrow}$}}
%                              \put(21,16){\vector(1,0){2}}
%                             \put(21,16){\line(1,0){2}}
%                             \put(21,16.5){\line(1,0){2}}
 \put(27,16){\makebox(0,0)[c]{$U_{T'}$}}
      \put(31,16){\makebox(0,0)[c]{${\dashrightarrow}$}}
%                              \put(30,16){\vector(1,0){2}}
 \put(35,16){\makebox(0,0)[c]{$U_{T''}$}}
                              \put(37,16){\vector(1,0){2}}

\end{picture}

Note that $\HT (V_{T''}[-1],U_T[1])=0$ as $V_{T''}\in \widecheck{\mr{add}{\omega}}$ and $U_T\in{_{\omega}\mc{X}}$, so there is a homomorphism $f: V_{T''}[-1]\to V_T$ such that the lower left square of the diagram commutates. Then we can easily complete it to be a commutative diagram of triangles for some $U_{T'}$ and $V_{T'}$. Since both $\widecheck{\mr{add}{\omega}}$ and $U_T\in{_{\omega}\mc{X}}$ are extension-closed, we have that $V_{T'}\in \widecheck{\mr{add}{\omega}}$ and $U_{T'}\in{_{\omega}\mc{X}}$, as desired.

(2) Let $X_T\to Y_T\to T\to$ and $X_{T''}\to Y_{T''}\to {T''}\to$ with $Y_T,Y_{T''}\in{\mc{X}_{\omega}}$ and $X_T,X_{T''}\in\widehat{\mr{add}{\omega}}$, given for $T, {T''}$ by Lemma \ref{Bl}(4). Consider the following diagram.

 \setlength{\unitlength}{0.09in}
 \begin{picture}(50,20)

                 \put(9,3.4){\vector(0,-1){2}}
                 \put(18,3.4){\vector(0,-1){2}}
                 \put(27,3.4){\vector(0,-1){2}}
                 \put(35,3.4){\vector(0,-1){2}}

 \put(9,5){\makebox(0,0)[c]{$T''[-1]$}}
                             \put(13,5){\vector(1,0){2}}
 \put(18,5){\makebox(0,0)[c]{$T$}}
                             \put(21,5){\vector(1,0){2}}
 \put(27,5){\makebox(0,0)[c]{$T'$}}
                             \put(30,5){\vector(1,0){2}}
 \put(35,5){\makebox(0,0)[c]{$T''$}}
                             \put(37,5){\vector(1,0){2}}

                 \put(9,9){\vector(0,-1){2}}
                 \put(18,9){\vector(0,-1){2}}
                 \put(27,9){\vector(0,-1){2}}
                 \put(35,9){\vector(0,-1){2}}
%                 \put(34.5,9){\line(0,-1){2}}
%                 \put(35,9){\line(0,-1){2}}

 \put(9,11){\makebox(0,0)[c]{$Y_{T''}[-1]$}}
      \put(14,11.5){\makebox(0,0)[c]{$\stackrel{f}{\dashrightarrow}$}}
%                             \put(13,11){\vector(1,0){2}}
%                             \put(14,12){\makebox(0,0)[c]{$f$}}
 \put(18,11){\makebox(0,0)[c]{$Y_T$}}
      \put(22,11){\makebox(0,0)[c]{${\dashrightarrow}$}}
%                             \put(21,11){\vector(1,0){2}}
 \put(27,11){\makebox(0,0)[c]{$Y_{T'}$}}
      \put(31,11){\makebox(0,0)[c]{${\dashrightarrow}$}}
%                             \put(30,11){\vector(1,0){2}}
 \put(35,11){\makebox(0,0)[c]{$Y_{T''}$}}
                             \put(37,11){\vector(1,0){2}}

                 \put(9,14.5){\vector(0,-1){2}}
                 \put(18,14.5){\vector(0,-1){2}}
                 \put(27,14.5){\vector(0,-1){2}}
                 \put(35,14.5){\vector(0,-1){2}}

 \put(9,16){\makebox(0,0)[c]{$X_{T''}[-1]$}}
      \put(14,16){\makebox(0,0)[c]{${\dashrightarrow}$}}
%                              \put(13,16){\vector(1,0){2}}
 \put(18,16){\makebox(0,0)[c]{$X_T$}}
      \put(22,16){\makebox(0,0)[c]{${\dashrightarrow}$}}
%                              \put(21,16){\vector(1,0){2}}
%                             \put(21,16){\line(1,0){2}}
%                             \put(21,16.5){\line(1,0){2}}
 \put(27,16){\makebox(0,0)[c]{$X_{T'}$}}
      \put(31,16){\makebox(0,0)[c]{${\dashrightarrow}$}}
%                              \put(30,16){\vector(1,0){2}}
 \put(35,16){\makebox(0,0)[c]{$X_{T''}$}}
                              \put(37,16){\vector(1,0){2}}

\end{picture}

Note that $\HT (Y_{T''}[-1],X_T[1])=0$ as $X_T\in \widehat{\mr{add}{\omega}}$ and $Y_{T''}\in{_{\omega}\mc{X}}$, so there is a homomorphism $f: Y_{T''}[-1]\to Y_T$ such that the lower left square of the diagram commutates. Then we can easily complete it to be a commutative diagram of triangles for some $X_{T'}$ and $Y_{T'}$. Since both $\widehat{\mr{add}{\omega}}$ and ${\mc{X}_{\omega}}$ are extension-closed, we have that $X_{T'}\in \widehat{\mr{add}{\omega}}$ and $Y_{T'}\in{_{\omega}\mc{X}}$, as desired.
\hfill$\Box$

%**************************************************
%%%%%%%%%%%%%%%%%%%%%%%%%%%%%%%%%%55

\vskip 10pt

Following the idea in [W1], we introduce the following notion in triangulated categories.

\bg{Def}\label{wG}%w-Goerenstin
Let $\mc{T}$ be a triangulated category and ${\omega}\subseteq\mc{T}$ be semi-selforthogonal. An object $G$ is called \textbf{$\omega$-Gorenstein} if $G\in {_{\omega}\mc{X}}\bigcap\mc{X}_{\omega}$.
\ed{Def}

%%%%%%%%%%%%%%%%%%%%%%%%%%%%%%%%%%%%%%%%%%%%%%%%%%%%%%%%%%%%%%%%%%%%%%%%%%%%%%%%%%

 The subcategory of all $\omega$-Gorenstein objects is denoted by $\mc{G}_{\omega}$ in the following. Hence, $\mc{G}_{\omega}={_{\omega}\mc{X}}\bigcap\mc{X}_{\omega}$. We note that  Asadollahi and Salarian [AS] already introduced the notion of $\xi$-Gorenstein objects in triangulated categories through relative homological algebra for triangulated algebras developed by Beligiannis [Be2].

 Let $\mc{C},\mc{D}$ be two subcategories of $\mc{T}$. A homomorphism $f: C\to D$ with $C\in\mc{C}$ and $D\in\mc{D}$ is said to be a {\it right $\mc{C}$-approximation} (or a $\mc{C}$-precover) of $D$ if $\HT (C',f)$ is surjective for any $C'\in\mc{C}$. Dually, a homomorphism $g: C\to D$ with $C\in\mc{C}$ and $D\in\mc{D}$ is said to be a {\it left $\mc{D}$-approximation} (or a $\mc{D}$-preenvelope) of $C$ if $\HT (g,D')$ is surjective for any $D'\in\mc{D}$.

%\vskip 5pt

We now recall the notion of $\mc{D}$-mutation pairs introduced by Iyama and Yoshino [IYo].
Let $\mc{U},\mc{V},\mc{D}$ be two subcategories of $\mc{T}$. Assume that $\mc{D}$ satisfies that $\HT (\mc{D},\mc{D}[1])=0$.  Following [IYo], we define $\mu^{-1}(\mc{U};\mc{D})$ to be the class of objects $T$ such that there is a triangle $U\stackrel{f}{\to} D\to T\to$ for some $U\in \mc{U}$ and $f$ is a left $\mc{D}$-approximation. Dually, we define $\mu(\mc{V};\mc{D})$ to be the class of objects $T$ such that there is a triangle $T\to D\stackrel{g}{\to} V\to $ for some $V\in \mc{V}$ and $g$ a right $\mc{D}$-approximation. The pair $(\mc{U},\mc{V})$ is a {\it{$\mc{D}$-mutation pair}} if $\mc{D}\subseteq \mc{U}\subseteq \mu(\mc{V};\mc{D})$ and $\mc{D}\subseteq \mc{V}\subseteq \mu^{-1}(\mc{U};\mc{D})$. Note that  $(\mc{U},\mc{V})$ is a $\mc{D}$-mutation pair if and only if $\mc{D}\subseteq \mc{U}\bigcap \mc{V}$ and both $\mc{U}=\mu(\mc{V};\mc{D})$ and $\mc{V}= \mu^{-1}(\mc{U};\mc{D})$ hold, by [IYo, Proposition 2.6].

Let $\mc{D}\subseteq \mc{C}$ be subcategories in $\mc{T}$. Following [IYo], we denote by $[\mc{D}]$ the ideal of homomorphisms in $\mc{C}$ which factors through objects in $\mc{D}$. The stable category of $\mc{C}$ by $\mc{D}$ is the quotient category $\mc{C}/[\mc{D}]$, for which objects are the same as $\mc{C}$, but homomorphisms are  homomorphisms in $\mc{C}$ modulo those in $[\mc{D}]$.

%%%%%%%%%%%%%%%%%%%%%%%%%%%%%%%%%%%%%%%%%%%%%%%%%%%%%%%%%%%%%%%%%%%%%%%%%%%%%%%%%%
%\hskip 18pt
%
\bg{Pro}\label{Dm}%D-mutation
Let ${\omega}$ be semi-selforthogonal. Let $\mc{D}=\mr{add}{\omega}$. Then

$(1)$ $\mc{G}_{\omega}$ is closed under extensions and direct summands and finite direct sums.

$(2)$ The pair $(\mc{G}_{\omega},\mc{G}_{\omega})$ is a $\mc{D}$-mutation pair.

$(3)$ The stable category $\underline{\mc{G}_{\omega}}:=\mc{G}_{\omega}/[\mc{D}]$ is a triangulated category.

$(4)$ For any integer $i$, $\underline{\mc{G}_{\omega[i]}}$ is naturally triangle equivalent to $\underline{\mc{G}_{\omega}}$.
\ed{Pro}

\Pf. (1) By Lemma \ref{Bl} (1).

(2) Take any $G\in\mc{G}_{\omega}$. Since $\mc{G}_{\omega}\subseteq \mc{X}_{\omega}$, there is a triangle $G\to M_0\stackrel{g}{\to} G_1\to$ with $M_0\in\mc{D}$ and $G_1\in\mc{X}_{\omega}$, by the definition. Note that $G\in {_{\omega}\mc{X}}\omega^{\perp}$, so it is easy to see that $g$ is a right $\mc{D}$-approximation. Moreover, since $M_0\in \mc{D}\subseteq  {_{\omega}\mc{X}}$ and ${_{\omega}\mc{X}}$ is coresolving by Lemma \ref{Bl} (2), we have that $G_1\in {_{\omega}\mc{X}}$ too. It follows that $G_1\in \mc{G}_{\omega}$. Thus, $\mc{G}_{\omega}\subseteq \mu(\mc{G}_{\omega};\mc{D})$ by the definition.

Dually, one can also prove that $\mc{G}_{\omega}\subseteq \mu^{-1}(\mc{G}_{\omega};\mc{D})$. Since $\mc{D}\subseteq\mc{G}_{\omega}$, we obtain that  $(\mc{G}_{\omega},\mc{G}_{\omega})$ is a $\mc{D}$-mutation pair.

(3) By [IYo, Section 4].

(4) Obviously.
\hfill $\Box$

%
%%%%%%%%%%%%%%%%%%%%%%%%%%%%%%%%%%%%%%%%%%%%%%%%%%%%%%%%%%%%%%%%%%%%%%%%%%%%%%%%%%%
%\hskip 18pt
%
%As shown in [IY], the shift functor $\langle 1\rangle$ of the triangulated category $\underline{\mc{G}_{\omega}}$ sends an object $G\in\underline{\mc{G}_{\omega}}$ to an object $\mu(G;\mc{D})$ (up to isomorphisms in $\underline{\mc{G}_{\omega}}$) and any triangle in $\underline{\mc{G}_{\omega}}$ can be obtain from a triangle in $\mc{T}$. We refer to the paper [IY] for other properties of $\mc{D}$-mutation pairs.
%
%**************************************************
%%%%%%%%%%%%%%%%%%%%%%%%%%%%%%%%%%55
\vskip 10pt

Recall that, for a triangulated subcategory $\mc{B}\subseteq \mc{T}$, the Verdier's quotient triangulated category $\mc{T}/\mc{B}$ is well-defined and there is a canonical functor $Q: \mc{T}\to \mc{T}/\mc{B}$ with the universal property that any functor $G$ from $\mc{T}$ to another category $\mc{A}$ with $G(\mc{B})=0$ actually factors through the functor $Q$.

 Following the idea in [C], we call $\mc{T}/\mc{B}$ the {\it relative singularity category} of $\mc{T}$ with respect to $\omega$, denoted by $\mc{T}_{\omega}$, in case $\mc{B}=\langle \mr{add}\omega\rangle$ for a semi-selforthogonal subcategory $\omega$ in $\mc{T}$.

%**************************************************
%%%%%%%%%%%%%%%%%%%%%%%%%%%%%%%%%%55
The following lemma  will be needed in the proof of our main theorem.

\bg{Lem}\label{Df}% D factor
Let ${\omega}\subseteq\mc{T}$ be semi-selforthogonal. Let $\mc{D}=\mr{add}{\omega}$. Then

$(1)$ Both $\HT (\widecheck{\mr{add}{\omega}}, {_{\omega}\mc{X}})$ and $\HT ({\mc{X}_{\omega}}, \widehat{\mr{add}{\omega}})$ are contained in $[\mc{D}]$.

$(2)$ $\mc{G}_{\omega}\bigcap\ \langle \mc{D}\rangle= \mc{D} $.
\ed{Lem}

\Pf. (1) Take any $f\in\HT (X,Y)$, where $X\in\widecheck{\mr{add}{\omega}}$ and $Y\in {_{\omega}\mc{X}}$. For $X\in \widecheck{\mr{add}{\omega}}$, there is a triangle $X\stackrel{h}\to M_X\to X'\to$ with $M_X\in\mr{add}{\omega}$ and $X'\in \widecheck{\mr{add}{\omega}}$. Since $\HT (X'[-1],Y)=0$ following from the fact that $\HT (\widecheck{\mr{add}{\omega}},{_{\omega}\mc{X}}[1])=0$, we have that $f$ factors through $M_X$ and hence $f\in [\mr{add}{\omega}]$.
The remained part can be proved dually.

(2) By [W2, Lemma 2.1 and Corollary 2.6].
\hfill $\Box$

%**************************************************
%%%%%%%%%%%%%%%%%%%%%%%%%%%%%%%%%%55

\vskip 10pt

We now present our main result as follows. Interesting applications will be given in next section.

\bg{Th}\label{T}%
Let $\mc{T}$ be a triangulated category and ${\omega}\subseteq\mc{T}$ be semi-selforthogonal. Let $\mc{D}=\mr{add}{\omega}$ and $\mc{G}_{\omega}= {_{\omega}\mc{X}}\bigcap\mc{X}_{\omega}$. Assume that $\widehat{\mc{X}_{\omega}}=\mc{T}=\widecheck{_{\omega}\mc{X}}$. Then there is a triangle equivalence between the stable category $\underline{\mc{G}_{\omega}}$ and the relative singularity category $\mc{T}_{\omega}$.
\ed{Th}

\Pf. We firstly construct a functor $F: \mc{T}_{\omega}\to \underline{\mc{G}_{\omega}}$ and then prove that it induces the desired triangle equivalence.

\vskip 5pt
{\bf Step 1}. Constructing the functor $F:\mc{T}_{\omega} \to\underline{\mc{G}_{\omega}}$.

Take any $T, T'\in\mc{T}$ and $f\in\HT (T,T')$. Since $\mc{T}=\widecheck{_{\omega}\mc{X}}$, there are triangles $U_T\to V_T\to T\to$ and $U_{T'}\to V_{T'}\to T'\to$ with $U_T,U_{T'}\in{_{\omega}\mc{X}}$ and $V_T,V_{T'}\in\widecheck{\mr{add}{\omega}}$, by Lemma \ref{Bl} (3). Note that $\HT (\widecheck{\mr{add}{\omega}}, {_{\omega}\mc{X}}[1])=0$, the homomorphism $f$ induces a homomorphism $f^V: V_T\to V_{T'}$ and consequently induces a homomorphism $f^U: U_T\to U_{T'}$, such that the following diagram is commutative.

 \setlength{\unitlength}{0.09in}
 \begin{picture}(50,8)

 \put(18,6){\makebox(0,0)[c]{$U_T$}}
                             \put(20,6){\vector(1,0){4}}
                                 \put(22,7){\makebox(0,0)[c]{$u$}}
 \put(27,6){\makebox(0,0)[c]{$V_T$}}
                             \put(29,6){\vector(1,0){4}}
                                 \put(31,7){\makebox(0,0)[c]{$v$}}
 \put(35,6){\makebox(0,0)[c]{$T$}}
                             \put(37,6){\vector(1,0){4}}
%\put(43,6){\makebox(0,0)[c]{$T$}}

                \put(18,5){\vector(0,-1){3}}
                  \put(19,4){\makebox(0,0)[c]{$f^U$}}
                \put(27,5){\vector(0,-1){3}}
                  \put(28,4){\makebox(0,0)[c]{$f^V$}}
                \put(35,5){\vector(0,-1){3}}
                  \put(36,4){\makebox(0,0)[c]{$f$}}
%                 \put(34.5,9){\line(0,-1){2}}
%                 \put(35,9){\line(0,-1){2}}

 \put(18,1){\makebox(0,0)[c]{$U_{T'}$}}
                             \put(20,1){\vector(1,0){4}}
                                 \put(22,2){\makebox(0,0)[c]{$u'$}}
 \put(27,1){\makebox(0,0)[c]{$V_{T'}$}}
                             \put(29,1){\vector(1,0){4}}
                                 \put(31,2){\makebox(0,0)[c]{$v'$}}
 \put(35,1){\makebox(0,0)[c]{$T'$}}
                             \put(37,1){\vector(1,0){4}}

\end{picture}
\vskip 5pt

Since $\widehat{\mc{X}_{\omega}}=\mc{T}$ too, for $U_T, U_{T'}$, there are triangles $X^U_T\to Y^U_T\to U_T\to$ and $X^U_{T'}\to Y^U_{T'}\to U_{T'}\to$ with $Y^U_T,Y^U_{T'}\in{\mc{X}_{\omega}}$ and $X^U_T,X^U_{T'}\in\widehat{\mr{add}{\omega}}$, by Lemma \ref{Bl} (4). Note also that $\HT ({\mc{X}_{\omega}},\widehat{\mr{add}{\omega}}[1])=0$, so there is the following commutative diagram of triangles.

 \setlength{\unitlength}{0.09in}
 \begin{picture}(50,8)

 \put(18,6){\makebox(0,0)[c]{$X^U_T$}}
                             \put(20,6){\vector(1,0){4}}
                                 \put(22,7){\makebox(0,0)[c]{$x$}}
 \put(27,6){\makebox(0,0)[c]{$Y^U_T$}}
                             \put(29,6){\vector(1,0){4}}
                                 \put(31,7){\makebox(0,0)[c]{$y$}}
 \put(35,6){\makebox(0,0)[c]{$U_T$}}
                             \put(37,6){\vector(1,0){4}}
%\put(43,6){\makebox(0,0)[c]{$T$}}

                \put(18,5){\vector(0,-1){3}}
                  \put(19,4){\makebox(0,0)[c]{$f^U_X$}}
                \put(27,5){\vector(0,-1){3}}
                  \put(28,4){\makebox(0,0)[c]{$f^U_Y$}}
                \put(35,5){\vector(0,-1){3}}
                  \put(36,4){\makebox(0,0)[c]{$f^U$}}
%                 \put(34.5,9){\line(0,-1){2}}
%                 \put(35,9){\line(0,-1){2}}

 \put(18,1){\makebox(0,0)[c]{$X^U_{T'}$}}
                             \put(20,1){\vector(1,0){4}}
                                 \put(22,2){\makebox(0,0)[c]{$x'$}}
 \put(27,1){\makebox(0,0)[c]{$Y^U_{T'}$}}
                             \put(29,1){\vector(1,0){4}}
                                 \put(31,2){\makebox(0,0)[c]{$y'$}}
 \put(35,1){\makebox(0,0)[c]{$U_{T'}$}}
                             \put(37,1){\vector(1,0){4}}

\end{picture}
\vskip 5pt

Since $\widehat{\mr{add}{\omega}}\subseteq {_{\omega}\mc{X}}$ and ${_{\omega}\mc{X}}$ is extension-closed, we see that $Y^U_T,Y^U_{T'}\in {_{\omega}\mc{X}}$ from the above triangles. It follows that $Y^U_T,Y^U_{T'}\in \mc{G}_{\omega}={_{\omega}\mc{X}}\bigcap \mc{X}_{\omega}$.

We now show that $S(f^U_Y)$ is unique only depending on the homomorphism $f$, where $S: \mc{G}_{\omega}\to \underline{\mc{G}_{\omega}}:=\mc{G}_{\omega}/[\mc{D}]$ is the canonical stabilize functor. Clearly, it is enough to prove that $S(f^U_Y)=0$ in $\underline{\mc{G}_{\omega}}$ whenever $f=0$, that is, it is enough to prove that $S(f^U_Y)\in [\mc{D}]$.

Note that we have the following commutative diagram of triangles.

 \setlength{\unitlength}{0.09in}
 \begin{picture}(50,8)

 \put(8,6){\makebox(0,0)[c]{$T[-1]$}}
                             \put(11,6){\vector(1,0){4}}
                                 \put(13,7){\makebox(0,0)[c]{$t$}}
 \put(18,6){\makebox(0,0)[c]{$U_T$}}
                             \put(20,6){\vector(1,0){4}}
                                 \put(22,7){\makebox(0,0)[c]{$u$}}
 \put(27,6){\makebox(0,0)[c]{$V_T$}}
                             \put(29,6){\vector(1,0){4}}
                                 \put(31,7){\makebox(0,0)[c]{$v$}}
 \put(35,6){\makebox(0,0)[c]{$T$}}
                             \put(37,6){\vector(1,0){4}}
%\put(43,6){\makebox(0,0)[c]{$T$}}

                \put(9,5){\vector(0,-1){3}}
                  \put(10,3.8){\makebox(0,0)[c]{$0$}}
                \put(18,5){\vector(0,-1){3}}
                  \put(19,4){\makebox(0,0)[c]{$f^U$}}
                \put(27,5){\vector(0,-1){3}}
                  \put(28,4){\makebox(0,0)[c]{$f^V$}}
                     \put(26,5){\vector(-2,-1){6}}
                     \put(23,4){\makebox(0,0)[c]{$\theta$}}
                \put(35,5){\vector(0,-1){3}}
                  \put(36,3.8){\makebox(0,0)[c]{$0$}}
%                 \put(34.5,9){\line(0,-1){2}}
%                 \put(35,9){\line(0,-1){2}}

 \put(8,1){\makebox(0,0)[c]{$T'[-1]$}}
                             \put(11,1){\vector(1,0){4}}
                                 \put(13,2){\makebox(0,0)[c]{$t'$}}
 \put(18,1){\makebox(0,0)[c]{$U_{T'}$}}
                             \put(20,1){\vector(1,0){4}}
                                 \put(22,2){\makebox(0,0)[c]{$u'$}}
 \put(27,1){\makebox(0,0)[c]{$V_{T'}$}}
                             \put(29,1){\vector(1,0){4}}
                                 \put(31,2){\makebox(0,0)[c]{$v'$}}
 \put(35,1){\makebox(0,0)[c]{$T'$}}
                             \put(37,1){\vector(1,0){4}}

\end{picture}
\vskip 5pt

Since $tf^U=0$, we obtain a homomorphism $\theta: V_T\to U_{T'}$ such that $f^U=u\theta$. Now we consider the following commutative diagram

 \setlength{\unitlength}{0.09in}
 \begin{picture}(50,12)

 \put(18,10){\makebox(0,0)[c]{$X^U_T$}}
                             \put(20,10){\vector(1,0){4}}
                                 \put(22,11){\makebox(0,0)[c]{$x$}}
 \put(27,10){\makebox(0,0)[c]{$Y^U_T$}}
                             \put(29,10){\vector(1,0){4}}
                                 \put(31,11){\makebox(0,0)[c]{$y$}}
 \put(35,10){\makebox(0,0)[c]{$U_T$}}
                             \put(37,10){\vector(1,0){4}}

                \put(18,9){\vector(0,-1){7}}
                  \put(19,6){\makebox(0,0)[c]{$f^U_X$}}
\put(27,9){\vector(0,-1){7}}
                  \put(28,6){\makebox(0,0)[c]{$f^U_Y$}}
                \put(35,9){\vector(0,-1){7}}
                  \put(36,6){\makebox(0,0)[c]{$f^U$}}
                  \put(36,9){\vector(1,-1){2.5}}
                    \put(38,8){\makebox(0,0)[c]{$u$}}
                  \put(40,6){\makebox(0,0)[c]{$V_T$}}
                    \put(39,5){\vector(-1,-1){3}}
                      \put(38,3){\makebox(0,0)[c]{$\theta$}}
                    \put(38.5,5){\vector(-4,-1){11}}
                      \put(32,4){\makebox(0,0)[c]{$\rho$}}

 \put(18,1){\makebox(0,0)[c]{$X^U_{T'}$}}
                             \put(20,1){\vector(1,0){4}}
                                 \put(22,2){\makebox(0,0)[c]{$x'$}}
 \put(27,1){\makebox(0,0)[c]{$Y^U_{T'}$}}
                             \put(29,1){\vector(1,0){4}}
                                 \put(31,2){\makebox(0,0)[c]{$y'$}}
 \put(35,1){\makebox(0,0)[c]{$U_{T'}$}}
                             \put(37,1){\vector(1,0){4}}
                                 \put(39,2){\makebox(0,0)[c]{$h$}}
 \put(44,1){\makebox(0,0)[c]{$X^U_{T'}[1]$}}
                             \put(46,1){\vector(1,0){4}}

\end{picture}
\vskip 5pt

Since $\HT (\widecheck{\mr{add}{\omega}},\widehat{\mr{add}{\omega}}[1])=0$, so $\theta h=0$ and then there is $\rho\in\HT (V_T, Y^U_{T'})$ such that $\rho y'=\theta$. It follows that $f^U_Yy'=yf^U=yu\theta=yu\rho y'$. Thus we have that $(f^U_Y-yu\rho)y'=0$. Hence there is $\delta: Y^U_T\to X^U_{T'}$ such that $f^U_Y-yu\rho=\delta x'$. It follows that $f^U_Y=-yu\rho +\delta x'$. But $\rho, \delta\in [\mc{D}]$ by Lemma \ref{Df}, it follows that $f^U_Y\in [\mc{D}]$ too.

Following from the above arguments, one can easily deduce that there is a well-defined functor $F: \mc{T}\to \underline{\mc{G}_{\omega}}$, by setting $F(T)=S(Y_T)$ for $T\in\mc{T}$ and $F(f)=S(f^U_Y)$ for $f\in\HT(T,T')$.

We claim that $F(T)=0$ for any $T\in\langle\mc{D}\rangle$ and hence $F$ factors through $\mc{T}_{\omega}$. The induced functor will still be denoted by $F$.

In fact, since both $\widecheck{\mr{add}{\omega}}$ and $\widehat{\mr{add}{\omega}}$ are contained in $\langle\mc{D}\rangle$ and  $\langle\mc{D}\rangle$ is a triangulated subcategory, it follows from the above construction of $F(T)$ that $Y^U_T\in \mc{G}_{\omega}\bigcap\ \langle\mc{D}\rangle$, in case  $T\in\langle\mc{D}\rangle$. But $\mc{G}_{\omega}\bigcap\ \langle\mc{D}\rangle=\mc{D}$ by Lemma \ref{Df}. So we have that $F(T)=S(Y^U_T)=0$ in $\underline{\mc{G}_{\omega}}$.
\vskip 5pt

{\bf Step 2}.  $F$ is a triangle functor.

Let $T\to T'\to T''\to$ be a triangle in $\mc{T}$. By Lemma \ref{Tl}, we have a triangle $Y^U_T\to Y^U_{T'}\to Y^U_{T''}\to$ in $\mc{G}_{\omega}$. It further induces a triangle $S(Y^U_T)\to S(Y^U_{T'})\to S(Y^U_{T''})\to$ in $\underline{\mc{G}_{\omega}}$ [IYo]. So we obtain that the functor $F$ induces a triangle $F(T)\to F(T')\to F(T'')\to$ by the definition.

\vskip 5pt

{\bf Step 3}.  $F$ is dense, full and faithful, so that $F$ induces an equivalence.

For any object $G\in\underline{\mc{G}_{\omega}}$, we can assume that $G\in\mc{G}_{\omega}\subseteq\mc{T}$ by the definition. Then it is easy to see that $F(G)=G$ from Step 1. So that $F$ is dense.

For any $T\in\mc{T}$, we see that $T\simeq Y^U_T$ in $\mc{T}_{\omega}$ from triangles $U_T\to V_T\to T\to$ and $X^U_T\to Y^U_T\to U_T\to$ in Step 1, since that $V_T\in\widecheck{\mr{add}{\omega}}$, $X^U_T\in\widehat{\mr{add}{\omega}}$ and that both $\widecheck{\mr{add}{\omega}}$ and $\widehat{\mr{add}{\omega}}$ are contained in $\langle\mc{D}\rangle$. It follows that $\mr{Hom}_{\mc{T}_{\omega}}(T,T')\simeq \mr{Hom}_{\mc{T}_{\omega}}(Y^U_T,Y^U_{T'})$, for any $T,T'\in \mc{T}_{\omega}$. Now take any $\underline{f}\in\mr{Hom}_{\underline{\mc{G}_{\omega}}}(F(T),F(T'))=\mr{Hom}_{\underline{\mc{G}_{\omega}}}(S(Y^U_T),S(Y^U_{T'}))$, then we see that there is a $f\in\HT(Y^U_T,Y^U_{T'})$ such that $S(f)=\underline{f}$ by the definition. Set $g=Q(f)$, where $Q: \mc{T}\to \mc{T}_{\omega}$ is the canonical functor. Then we may assume that $g\in \mr{Hom}_{\mc{T}_{\omega}}(T,T')$ by the above arguments. Note that $F(g)=\underline{f}$, so we have that $F$ is full.

To prove that $F$ is faithful, it is enough to show that $F$ is faithful on objects by [KZ, Page 45], since $F$ is a full triangle functor. Now let $F(T)=0$ in $\underline{\mc{G}_{\omega}}$. Then we have that $Y^U_T\in\mc{D}$ by the definition of $F$ in Step 1. But it indicates that $T\in\langle\mc{D}\rangle$ from triangles in Step 1. It follows that $T\simeq 0$ in $\mc{T}_{\omega}$. Hence, $F$ is faithful.
\hfill$\Box$
%
%
%%%%%%%%%%%%%%%%%%%%%%%%%%%%%%%%%%55---------------------------------------------------------------------------------------------------------------

\vskip 10pt

We thank Dong Yang for communicating us the following result [IYa,  Theorem 4.7] which was recently obtained by Iyama and Yang. Note that the notion of presilting subcategories in their paper is just the notion of semi-selforthogonal subcategories in our sense.

\bg{Cor}\label{Cy}% Cor yang dong
Let $\omega\subseteq \mc{T}$ be semi-selforthogonal. Assume that (P1) $\omega$ is functorially finite in $\mc{T}$ and (P2) for any $X\in\mc{T}$, $\HT (X,\omega[i])=0=\HT (\omega, X[i])$ for $i\gg 0$. Then there is a triangle equivalence between $\underline{^{\perp}\omega\bigcap\omega^{\perp}}$ and relative singularity category $\mc{T}_{\omega}$.
\ed{Cor}

\Pf. We will prove that $(P1)$ implies $\omega^{\perp}={_{\omega}\mc{X}}$ and ${^{\perp}\omega}={\mc{X}_{\omega}}$, and that together with (P2) implies that $\mc{T}=\langle{_{\omega}\mc{X}}\rangle=\langle\mc{X}_{\omega}\rangle$. So the conclusion follows from Theorem \ref{T}.

Indeed, take any $X\in\omega^{\perp}$ and set $X_0:=X$. Since $\omega$ is functorially finite, there is a homomorphism $f_0: M_0\to X$ with $M_0\in\mr{add}\omega$ such that $\HT (\omega, f_0)$ is surjective. So $f_0$ induces a triangle $X_1\to M_0\to X_0\to$ with $X_1\in\omega^{\perp}$. Repeating the process to $X_1$, and so on, we see that $X=X_0\in {_{\omega}\mc{X}}$ actually by the definition. Dually, one also has that ${^{\perp}\omega}={\mc{X}_{\omega}}$.

Take any $X\in\mc{T}$, then there is some integer $n_X$ such that $\HT (X,\omega[i])=0$, for all $i>n_X$, by the condition (P2). If $n_X\le 0$, then we obtain that $X\in {^{\perp}\omega}$. If $n_X>0$, then we obtain that $X[-n_X]\in{^{\perp}\omega}$. In both cases, we have that $X\in\langle{^{\perp}\omega}\rangle=\langle\mc{X}_{\omega}\rangle$ by the above arguments. It follows that $\mc{T}=\langle\mc{X}_{\omega}\rangle$. Similarly, one can also obtain that $\mc{T}=\langle{_{\omega}\mc{X}}\rangle$.
\hfill$\Box$
%%%%%%%%%%%%%%%%%%%%%%%%%%%%%%%%%%55
%%%%%%%%%%%%%%%%%%%%%%%%%%%%%%%%%%55--------------------------------------------------------------------------------------------------------------------
%%%%%%%%%%%%%%%%%%%%%%%%%%%%%%%%%%55
\bg{Rem}\label{Rm} %remark
%$(1)$
Note that for a ring $R$, the semi-selorthogonal subcategory $\mr{Proj}=\mr{Add}R$ is usually not covariantly finite (hence not funtorially finite) in $\mr{Mod}R$, unless it is product-complete which means that any product of projective modules is a direct summand of some direct sum of projective (and hence also projective) [RS]. It only happens when $R$ is right perfect.
%$(2)$ Our corollaries
%
\ed{Rem}

\bskip\
\vskip 20pt
% ----------------------------------------------------------------------
%% ----------------------------------------------------------------------
%\def\baselinestretch{1}

\section{Applications}
%%
%%
%%%% ----------------------------------------------------------------------
%%

\vskip 2pt
%\mskip\
\hskip 10pt

Let $R$ be a ring. As usual, we denote by $\mr{Mod}R$ ($\mr{mod}R$, $\mr{Proj}$, $\mr{Inj}$, $\mc{P}$, resp.) the category of all (finitely presented, projective, injective, finitely generated projective, resp.) $R$-modules. And we denote by $\mc{D}^b(\mr{Mod}R)$ the derived category of complexes over $\mr{Mod}R$  and by $\mc{K}^b(\mr{Proj})$ ($\mc{K}^b(\mc{P})$, resp.) the homotopy category of complexes over $\mr{Proj}$ ($\mc{P}$, resp.).  If $R$ is left coherent, then the category $\mr{mod}R$ is an abelian category and we denote by $\mc{D}^b(R)$ the derived category of complexes over $\mr{mod}R$. By identifying an $R$-module $M$ with the complex with terms $M$ in 0-th position and 0 otherwise, we view $\mr{Mod}R$ as a full subcategory of $\mc{D}^b(\mr{Mod}R)$. For a complex $M\in \mc{D}^b(\mr{Mod}R)$, we denote by $\mr{Add}M$ the subcategory of all direct summands of direct sums of $M$. Thus, $\mr{Proj}=\mr{Add}R$. Clearly, we also have that $\mc{P}=\mr{add}R$.

Recall from [EJ] that an $R$-module $X$ is (finitely generated) {\it Gorenstein projective} provided that there is an infinite exact sequence $\cdots\to P_{-1}\to P_0\stackrel{f}{\to} P_1\to\cdots\ $ of (finitely generated) projective modules such that $X\simeq \mr{Ker}f$ and the sequence induced by the functor $\mr{Hom}_R(-,P)$ is also exact for any projective module $P$. Dually, an $R$-module $X$ is {\it Gorenstein injective} provided that there is an infinite exact sequence $\cdots\to I_{-1}\to I_0\stackrel{f}{\to} I_1\to\cdots\ $ of injective modules such that $X\simeq \mr{Ker}f$ and the sequence induced by the functor $\mr{Hom}_R(I,-)$ is also exact for any injective module $I$. The subcategories of  (finitely generated) Gorenstein projective and Gorenstein injective modules, respectively, are denoted by $\mc{G}P$ ($\mc{G}p$) and $\mc{G}I$ respectively.

Note that every module has a Gorenstein projective resolution and a Gorenstein injective resolution and that Gorenstein projective dimensions of modules are well-defined, as well as Gorenstein injective dimensions of modules. Moreover, the supersum of all Gorenstein  projective dimensions of modules coincides with the supersum of all Gorenstein injective dimensions of modules, which is called the {\it Gorenstein global dimension} of the ring [Be1,BM].

An important class of rings with finite Gorenstein global dimension are Gorenstein rings, that is, left and right noetherian rings with finite left and right self-injective dimensions. However, there are also non-noetherian rings with finite Gorenstein global dimension, as shown by the following example in [MT], where the ring $R\bowtie E$ is defined for a ideal $E$ of the ring $R$ to be a subring of $R\times R$ with the   multiplicative
structure given by $(r, e)(s, f):=(rs, rf+se+ef)$, where $r, s \in R$ and $e, f\in E$.

%%%%%%%%%%%%%%%%%%%%%%%%%%%%%%%%%%55
%\vskip 10pt

\bg{Exm}\label{Eg}% Example Gorenstein gd finite {\mr [MT]}
Let $n>0$ be an integer and let $R_n=R[X_1,X_2,\cdots\ , X_n]$ be the
polynomial ring in $n$ indeterminates over a commutative non-noetherian hereditary ring
$R$. Let $S_n:= R_n\bowtie R_nX_1$. Then $S_n$ is a non-noetherian coherent ring with Gorenstein global dimension $n+1$.

\ed{Exm}

We note that rings of finite Gorenstein global dimension were called Gorenstein rings in [Be1]. To distinguish with the usual Gorenstein rings (i.e., rings with two-sided noetherian and two-sided finite self-injective dimension), we reserve the name of rings of Gorenstein global dimension.

Over a ring of finite Gorenstein global dimension, we have the following characterizations of Gorenstein projective modules and Gorenstein injective modules, see [EJ, BM].

%%%%%%%%%%%%%%%%%%%%%%%%%%%%%%%%%%55
%\vskip 10pt

\bg{Lem}\label{Gl}% Gorenstein lemma
Let $R$ be a ring of finite Gorenstein global dimension. Then $\mc{G}P=\mr{KerExt}^{i>0}_R(-,\mr{Proj})$ and $\mc{G}I=\mr{KerExt}^{i>0}_R(\mr{Inj},-)$. Moreover, if $R$ is left coherent, then it also holds that $\mc{G}p=\mc{G}P\bigcap \mr{mod}R=\mr{KerExt}^{i>0}_R(-,R)\bigcap \mr{mod}R$.
\ed{Lem}

%%%%%%%%%%%%%%%%%%%%%%%%%%%%%%%%%%55---------------------------------------------------------------------------------------------------------------
%\vskip 10pt
We have the following result which describes the subcategory $\mc{G}_{\omega}$ in derived categories over rings of finite Gorenstein global dimension, in case $\omega=\mr{Proj}$, $\mr{Inj}$, $\mc{P}$, respectively.

\bg{Lem}\label{Gfl}% Gorenstein finite lemma
Let $R$ be a ring of finite Gorenstein global dimension. Then  $\mc{G}P=\mc{G}_{\mr{Proj}}$ and $\mc{G}I=\mc{G}_{\mr{Inj}}$, considered in $\mc{D}^b(\mr{Mod}R)$.
Moreover, if $R$ is also a left coherent ring, then it also holds that $\mc{G}p=\mc{G}_{R}$, considered in $\mc{D}^b(R)$.
\ed{Lem}

\Pf. We prove the last part. The first and second parts can be proved in a similar way and in a dual way, respectively.

 The part ${\mc{G}p}\subseteq {\mc{G}_R}$ is clear. On the other hand, take a complex $X\in {\mc{G}_R}$. Then $X\in {_R\mc{X}}\subseteq R^{\perp}$ and $X\in\mc{D}^b(R)\simeq \mc{K}^{-,b}(\mc{P})$, so there is an $R$-module $X_n$ and finitely generated projective modules $P_i$, $0\le i\le n$, such that there are triangles $X_i\stackrel{f_i}{\to} P_i\to X_{i-1}\to$ for $1\le i\le n$, where $X_0:=X$. Since $X\in\mc{G}_R\subseteq {^{\perp}R}$ too, we easily obtain that each $X_i\in {^{\perp}R}$. In particular, we have that $X_n\in {^{\perp}R}\bigcap\mr{mod}R=\mr{KerExt}^{i>0}_R(-,R)\bigcap\mr{mod}R=\mc{G}p$. It follows that there is a monomorphism $f: X_n\to Q$ for some projective $R$-module $Q$, by the definition. But $\mr{Hom}_R(f_n,Q)\simeq\HD(f_n,Q)$ is epic, so that $f$ factors through $f_n$. Then $f_n\in \HD(X_n,Q)\simeq\mr{Hom}_R(X_n,Q)$ is also a monomorphism of $R$-modules. Thus we have an exact sequence $0\to X_n\stackrel{f_n}{\to} P_n\to \mr{Cok}f_n\to 0$ which further gives a triangle $X_n\stackrel{f_n}{\to} P_n\to \mr{Cok}f_n\to$. Since $X_n\stackrel{f_n}{\to} P_n\to X_{n-1}\to$ is also a triangle, we obtain that $X_{n-1}\simeq \mr{Cok}f_n$ is an $R$-module. Repeating the process to $X_{n-1}$, and so on, we finally obtain that $X\in\mr{mod}R$. It follows that $X\in {^{\perp}R}\bigcap \mr{mod}R = \mr{KerExt}^{i>0}_R(-,R)\bigcap \mr{mod}R=\mc{G}p$.
\hfill$\Box$

%%%%%%%%%%%%%%%%%%%%%%%%%%%%%%%%%%55---------------------------------------------------------------------------------------------------------------
\vskip 10pt

Recall that the big singularity category $\mc{D}_{\mr{Sg}}(R)=\mc{D}^b(\mr{Mod}R)/\mc{K}^b(\mr{Proj})$ and the singularity category $\mc{D}_{sg}(R)=\mc{D}^b(R)/\mc{K}^b(\mc{P})$. Note that $\mc{K}^b(\mr{Proj})=\langle \mr{Proj}\rangle$ and $\mc{K}^b(\mc{P})=\langle \mc{P}\rangle=\langle R\rangle$.

As a corollary of Theorem \ref{T}, we have the following results. Note that (1) is a theorem of Beligiannis [Be1] and that (2) extends the classical result to coherent rings, see [Be1, Bu, C, R2] etc. %Our settings seem more general.

\bg{Cor}\label{Cfg}% Cor finite Gorenstein
$(1)$ Let $R$ be a ring of finite Gorenstein global dimension. Then there are triangle equivalences between $\underline{\mc{G}P}$, $\underline{\mc{G}I}$ and the big singularity category $\mc{D}_{\mr{Sg}}(R)$.

$(2)$ Let $R$ be a left coherent ring of finite Gorenstein global dimension. Then there is a triangle equivalence between $\underline{\mc{G}p}$ and the singularity category $\mc{D}_{sg}(R)$.
\ed{Cor}

\Pf. (1)  Since $R$ has finite Gorenstein global dimension, the class of modules of finite projective dimension coincides with the class of modules of finite injective dimension, see [BM]. So that $\langle \mr{Proj}\rangle = \langle \mr{Inj}\rangle$ in $\mc{D}^b(\mr{Mod}R)$.

Let $\omega=\mr{Proj}$ and $\mc{T}=\mc{D}^b(\mr{Mod}R)$. Then $\mc{T}_{\omega}=\mc{D}_{\mr{Sg}}(R)$.
By Lemma \ref{Gfl}, we have that ${\mc{G}_{\omega}}={\mc{G}P}$. Since every $R$-module has finite Gorenstein projective dimension, we obtain that $\mr{Mod}R\subseteq \langle\mc{G}P\rangle$. So we have that $\mc{D}^b(\mr{Mod}R)=\langle\mc{G}P\rangle$. Note that ${\mc{G}P}={\mc{G}_{\omega}}={_{\omega}\mc{X}}\bigcap\mc{X}_{\omega}\subseteq \mc{D}^b(\mr{Mod}R)$,  these imply that $\mc{D}^b(\mr{Mod}R)=\langle{_{\omega}\mc{X}}\rangle=\langle\mc{X}_{\omega}\rangle$. Combining with Lemma \ref{Rl}, we see that Theorem \ref{T} applies and then we obtain that $\underline{\mc{G}P}$ is triangle equivalent to the big singularity category $\mc{D}_{\mr{Sg}}(R)$.

Similarly, we have that $\underline{\mc{G}I}$ is also triangle equivalent to the big singularity category $\mc{D}_{\mr{Sg}}(R)$.

(2) The proof of (1) restricts to $\mc{D}^b(R)$ well in the settings.
\hfill$\Box$

%%%%%%%%%%%%%%%%%%%%%%%%%%%%%%%%%%55---------------------------------------------------------------------------------------------------------------

\vskip 10pt

Let $\mc{T}$ be a triangulated category and $\omega\subseteq \mc{T}$. Recall that $\omega$ is {\it silting} if $\omega$ is semi-selforthogonal and $\langle \mr{add}\omega\rangle=\mc{T}$ [AI]. Let $R$ be a ring. A complex $M$ in $\mc{D}^b(R)$ ($\mc{D}^b(\mr{Mod}R)$, resp.) is called silting (big silting, resp.) if $\omega=\mr{add}M$ ($\omega=\mr{Add}M$, resp.) is silting in $\mc{K}^b(\mc{P})$ ($\mc{K}^b(\mr{Proj})$, resp.) [W2]. Clearly, silting complexes are always big silting.

We note that, up to shifts, a complex $M$ over a ring $R$ is big silting if and only if (1) $M\in\widehat{\mr{Add}R}$, (2) $M$ is coproduct-semi-selforthogonal, and (3) $R\in\widecheck{\mr{Add}M}$ [W2].

The following lemma is useful.
%%%%%%%%%%%%%%%%%%%%%%%%%%%%%%%%%%55---------------------------------------------------------------------------------------------------------------

\bg{Lem}\label{Sl}% silting lemma
Let $\mc{T}$ be a triangulated category and let $\mc{U},\mc{V}\subseteq \mc{T}$ be both semi-selforthogonal. If $\mc{U}\subseteq \mc{X}_{\mc{V}}$ and ${^{\perp}\mc{U}}\subseteq{^{\perp}\mc{V}}$, then $\mc{X}_{\mc{U}}\subseteq\mc{X}_{\mc{V}}$.
\ed{Lem}

\Pf. Take any $T\in \mc{X}_{\mc{U}}$, then there are triangles   $T_i\to U_i\to T_{i+1}\to$ for all $i\ge 0$, such that $T_0:=T$, $U_i\in \mr{add}\mc{U}$ and $T_i\in {^{\perp}{\mc{U}}}$, for each $i\ge 0$. Note that $\mc{U}\subseteq \mc{X}_{\mc{V}}$, so we have all $U_i\in \mr{add}\mc{U}\subseteq \mc{X}_{\mc{V}}$ by Lemma \ref{Bl} (1). Thus, for each $i$, we have triangles   $\natural_j: U^j_i\to V^j_i\to U^{j+1}_{i}\to$ for all $j\ge 0$, such that $U^0_i:=U_i$, $V^j_i\in \mr{add}\mc{V}$ and $U^j_i\in {^{\perp}{\mc{V}}}$, for each $j\ge 0$.

Consider firstly $i=0$. Note that we have the following induced commutative diagram of triangles, for some $T^1_0$.

\mskip\

 \setlength{\unitlength}{0.09in}
 \begin{picture}(50,18)

%                 \put(18,3.4){\vector(0,-1){2}}
                 \put(27,3.4){\vector(0,-1){2}}
                 \put(35,3.4){\vector(0,-1){2}}

% \put(18,5){\makebox(0,0)[c]{$U$}}
%                             \put(21,5){\vector(1,0){2}}
 \put(27,5){\makebox(0,0)[c]{$U^1_0$}}
                             \put(30,4.9){\line(1,0){2}}
                             \put(30,5.2){\line(1,0){2}}
 \put(35,5){\makebox(0,0)[c]{$U^1_0$}}

%                 \put(18,9){\vector(0,-1){2}}
                 \put(27,9){\vector(0,-1){2}}
                 \put(35,9){\vector(0,-1){2}}
%                 \put(34.5,9){\line(0,-1){2}}
%                 \put(35,9){\line(0,-1){2}}

 \put(18,11){\makebox(0,0)[c]{$T_0$}}
                             \put(21,11){\vector(1,0){2}}
 \put(27,11){\makebox(0,0)[c]{$V^0_0$}}
                             \put(30,11){\vector(1,0){2}}
 \put(35,11){\makebox(0,0)[c]{$T^1_0$}}
                             \put(37,11){\vector(1,0){2}}

                 \put(17.9,14.5){\line(0,-1){2}}
                 \put(18.2,14.5){\line(0,-1){2}}
                 \put(27,14.5){\vector(0,-1){2}}
                 \put(35,14.5){\vector(0,-1){2}}

 \put(18,16){\makebox(0,0)[c]{$T_0$}}
                              \put(21,16){\vector(1,0){2}}
%                             \put(21,16){\line(1,0){2}}
%                             \put(21,16.5){\line(1,0){2}}
 \put(27,16){\makebox(0,0)[c]{$U_0$}}
                              \put(30,16){\vector(1,0){2}}
 \put(35,16){\makebox(0,0)[c]{$T_1$}}
                              \put(37,16){\vector(1,0){2}}

\end{picture}
\sskip\

Denote by $\nabla_0$ the triangle $ T_0\to V^0_0\to T^1_0\to$ in the middle row and by $\sharp_1$  the triangle $ T_1\to T^1_0\to U^1_0\to$ in the right column in the above diagram. Since $T_1\in {^{\perp}{\mc{U}}}\subseteq {^{\perp}{\mc{V}}}$ and $U^1_0\in {^{\perp}{\mc{V}}}$, we obtain that $T^1_0\in {^{\perp}{\mc{V}}}$.
Similarly as the triangle $\nabla_0$ for $T_0$, we see that, for $T_1$, there is also a triangle $\dag_1:   T_1\to V^0_1\to T^1_1\to$ with $T^1_1\in {^{\perp}{\mc{V}}} $.
%
%Similarly, we see that, for each $i$, there is a triangle $T_i\to V^0_i\to T^1_i\to$ with $T^1_i\in {^{\perp}{\mc{V}}}$.
%

Three triangles $\dag_1, \sharp_1, \natural_1$ induce the following commutative diagram of triangles with $V^1_{T}=V^0_1\oplus V^1_0$, for some $T^2_0$, since $U^1_0\in {^{\perp}{\mc{V}}}$.

\mskip\

 \setlength{\unitlength}{0.09in}
 \begin{picture}(50,18)

                 \put(18,3.4){\vector(0,-1){2}}
                 \put(27,3.4){\vector(0,-1){2}}
                 \put(35,3.4){\vector(0,-1){2}}

 \put(18,5){\makebox(0,0)[c]{$U^1_0$}}
                             \put(21,5){\vector(1,0){2}}
 \put(27,5){\makebox(0,0)[c]{$V^1_0$}}
                             \put(30,5){\vector(1,0){2}}
 \put(35,5){\makebox(0,0)[c]{$U^2_0$}}
                             \put(37,5){\vector(1,0){2}}

                 \put(18,9){\vector(0,-1){2}}
                 \put(27,9){\vector(0,-1){2}}
                 \put(35,9){\vector(0,-1){2}}
%                 \put(34.5,9){\line(0,-1){2}}
%                 \put(35,9){\line(0,-1){2}}

 \put(18,11){\makebox(0,0)[c]{$T^1_0$}}
                             \put(21,11){\vector(1,0){2}}
 \put(27,11){\makebox(0,0)[c]{$V^1_{T}$}}
                             \put(30,11){\vector(1,0){2}}
 \put(35,11){\makebox(0,0)[c]{$T^2_0$}}
                             \put(37,11){\vector(1,0){2}}

                 \put(18,14.5){\vector(0,-1){2}}
                 \put(27,14.5){\vector(0,-1){2}}
                 \put(35,14.5){\vector(0,-1){2}}

 \put(18,16){\makebox(0,0)[c]{$T_1$}}
                              \put(21,16){\vector(1,0){2}}
%                             \put(21,16){\line(1,0){2}}
%                             \put(21,16.5){\line(1,0){2}}
 \put(27,16){\makebox(0,0)[c]{$V^0_1$}}
                              \put(30,16){\vector(1,0){2}}
 \put(35,16){\makebox(0,0)[c]{$T^1_1$}}
                              \put(37,16){\vector(1,0){2}}

\end{picture}

Thus, we obtain triangles $\nabla_1: T^1_0\to V^1_{T}\to T^2_0\to $ for $T^1_0$ and $\sharp_2: T^1_1\to T^2_0\to U^2_0\to $ for $T^1_1$.
It is easy to see that $T^2_0\in  {^{\perp}{\mc{V}}}$, since $T^1_1, U^2_0\in {^{\perp}{\mc{V}}}$.
%Similarly, for each $i$, there is a triangle $T^1_i\to V^1_{T_i}\to T^2_i\to $ with $V^1_{T_i}\in \mr{add}\mc{V}$ and $T^2_i\in {^{\perp}{\mc{V}}}$.
Similarly as $\nabla_1$ for $T^1_0$, there is a triangle $\dag_2: T^1_1\to V^1_{T'}\to T^2_1\to $ for some $V^1_{T'}\in \mr{add}\mc{V}$ and $T^2_1\in {^{\perp}{\mc{V}}}$, for  $T^1_1$.

Repeating the above progress to triangles $\dag_2, \sharp_2, \natural_2$, and so on, we finally obtain triangles $\nabla_j: T^j_0\to V^j_T\to T^{j+1}_0\to$ for all $j\ge 0$, where $T^0_0:=T$, $V^0_T:=V^0_0$, such that $V^j_T\in\mr{add}\mc{V}$ and $T^j_0\in{^{\perp}{\mc{V}}}$, for each $j>0$. Hence, we conclude that $T\in \mc{X}_{\mc{V}}$ by the definition.
\hfill$\Box$
%%%%%%%%%%%%%%%%%%%%%%%%%%%%%%%%%%55---------------------------------------------------------------------------------------------------------------

\bg{Cor}\label{Sc}% silting c
$(1)$ Let $R$ be a ring of finite Gorenstein global dimension and let $M$ be a big silting complex in $\mc{D}^b(\mr{Mod}R)$. Then there is a triangle equivalence between $\underline{\mc{G}_{\mr{Add}M}}$ and the big singularity category $\mc{D}_{\mr{Sg}}(R)$.

$(2)$ Let $R$ be a left coherent ring of finite Gorenstein global dimension and let $M$ be a silting complex in $\mc{D}^b(R)$. Then there is a triangle equivalence between $\underline{\mc{G}_M}$ and the singularity category $\mc{D}_{sg}(R)$.
\ed{Cor}

\Pf. (1) Let $\mc{T}:=\mc{D}^b(\mr{Mod}R)$ and $\omega:=\mr{Add}M$. Since $M$ is big silting, $\langle \omega\rangle=\langle \mr{Add}M\rangle=\mc{K}^b(\mr{Proj})$. Hence, $\mc{T}_{\omega}=\mc{D}_{\mr{Sg}}(R)$ by the involved definitions. By Theorem \ref{T} and Lemma \ref{Rl}, we need only to show that $\mc{D}^b(\mr{Mod}R)=\langle\mc{X}_{\omega} \rangle=\langle{_{\omega}\mc{X}}\rangle$.

Take any $X\in\mc{D}^b(\mr{Mod}R)\simeq \mc{K}^{-,b}(\mr{Proj})$. Since $M\in\mc{K}^b(\mr{Proj})$, there is some integer $n>0$ such that $X[n]\in M^{\perp}$. But $M^{\perp}={_{\mr{Add}M}\mc{X}}$ by [W2, Lemma 3.11]. It follows that $X\in\widecheck{{_{\mr{Add}M}\mc{X}}}$ from the fact $0\in {_{\mr{Add}M}\mc{X}}$ and triangles $X[i]\to 0\to X[i+1]\to $ for all $i\ge 0$. So we easily obtain that $\mc{D}^b(\mr{Mod}R)=\widecheck{{_{\mr{Add}M}\mc{X}}}=\langle{_{\omega}\mc{X}}\rangle$.

By Proposition \ref{Dm} (4), we may assume that $M\in R^{\perp}$. In this case, we have that $M\in\widehat{\mr{Add}R}$ and $R\in\widecheck{\mr{Add}M}$, since $M$ is silting [W2, Theorem 3.5]. It follows that ${^{\perp}(\mr{Proj})}={^{\perp}(\mr{Add}R)}\subseteq {^{\perp}(\mr{Add}M)}$ and $\mr{Proj}=\mr{Add}R\subseteq \widecheck{\mr{Add}M}\subseteq \mc{X}_{\mr{Add}M}$. Hence, Lemma \ref{Sl} applies and we have that $\mc{X}_{\mr{Proj}}\subseteq \mc{X}_{\mr{Add}M}$.  Since $R$ has finite Gorenstein global dimension, we have that $\mc{D}^b(\mr{Mod}R)=\langle \mc{G}P\rangle$ as the arguments in the proof of Corollary \ref{Cfg} (1). The above arguments together with the fact $\mc{G}P\subseteq \mc{X}_{\mr{Proj}}$ implies that $\mc{D}^b(\mr{Mod}R)=\langle \mc{G}P\rangle \subseteq \langle \mc{X}_{\mr{Proj}}\rangle \subseteq \langle \mc{X}_{\mr{Add}M}\rangle \subseteq \mc{D}^b(\mr{Mod}R)$. It follows that $\mc{D}^b(\mr{Mod}R)=\langle \mc{X}_{\mr{Add}M}\rangle$.

(2) The arguments in the proof of (1) can be well restricted to $\mc{D}^b(R)$ in case $R$ is left coherent.
\hfill$\Box$
%
%%%%%%%%%%%%%%%%%%%%%%%%%%%%%%%%%%%%55

\vskip 10pt

We end the paper with the following example.

\vskip 10pt

\noindent {\bf Example 3.9}\hskip 10pt Let $R$ be the path algebra given by the quiver \hskip 5pt $1^{_{\ \longrightarrow}}_{^{\ \longleftarrow}}\ 2$ \hskip 5pt with radical square zero. Then $R$ is a self-injective algebra without nontrivial tilting complexes. By [R2], the singularity category $\mc{D}_{sg}(R)$ is triangle equivalent to the stable category  $\underline{\mr{mod}}R$. The indecomposable objects in $\underline{\mr{mod}}R$ are $\{1,2\}$ (up to isomorphisms). Let $M$ be the  two-terms silting complex {\footnotesize $\begin{array}{c}   2\\ 1 \end{array} \to \begin{array}{c}   1\\ 2 \end{array}$} with {\footnotesize $\begin{array}{c}   1\\ 2 \end{array}$} in the 0-th position. The indecomposable objects in $\underline{\mc{G}_M}$ are $\{1, 1[1]\}$ (up to isomorphisms). By Theorem \ref{T}, there is a triangle equivalence between $\mc{D}_{sg}(R)$ and $\underline{\mc{G}_M}$. Indeed, we have that  $F(1)=1$ and $F(2)=1[1]$.

\vskip 30pt

\noindent {\bf Reference}
{\small
%
%\begin{thebibliography}{17}

%
   \begin{itemize}

\item[{[AI]}]  T. Aihara and O. Iyama, Silting mutation in triangulated categories, J. Lon. Math. Soc. 85 (3) (2012), 633-668.
\item[{[AS]}]  A. Asadollahi and S. Salarian, Gorenstein objects in triangulated categories, J. Algebra 281 (2004), 264-286.
\item[{[AB]}]  M. Auslander and M. Bridger, Stable module theory, Mem. Amer. Math. Soc. 94 (1969).
\item[{[AR]}] M. Auslander and  I. Reiten, Applications of contravariantly
finite subcategories, Adv. Math. 86 (1991), 111-152.
\item[{[AF]}]  L. Avramov and H. Foxby, Ring homomorphisms and finite Gorenstein dimension,  Proc. London
Math. Soc. 75 (3) (1997), 241-270.
\item[{[Be1]}] A. Beligiannis, The homological theory of contravariantly finite subcategories: Auslander-Buchweitz contexts, Gorenstein categories and (co)stabilization. Comm. Algebra 28 (10) (2000), 4547-4596.
\item[{[Be2]}] A. Beligiannis, Relative homological algebra and purity in triangulated categories, J. Algebra 227 (1) (2000),
268-361.
\item[{[BM]}] D. Bennis and N. Mahdou, Global Gorenstein dimensions, Proc. Amer. Math.
Soc. 138 (2) (2010), 461-465.
\item[{[Bu]}] R. Buchweitz, Maximal Cohen-Macaulay modules and Tate cohomology over Gorenstein rings. Hamburg, p. 155 (1987) (unpublished manuscript).
\item[{[C]}] X. Chen, Relative singularity categories and Gorenstein projective modules, Math Nachr. 284 (2-3) (2011), 19-212.
\item[{[CZ]}] X. Chen and P. Zhang, Quotient triangulated categories, Manuscripta Math. 123 (2007), 167-183.
\item[{[EJ]}] E. Enochs and O. Jenda, Relative homological algebra, Walter de Gruyter, Berlin-New York, 2000.
\item[{[H1]}]   D. Happel, Triangulated Categories in the Representation
Theory of Finite Dimensional Algebras,  London Math. Soc.
Lect. Note Ser.  119  (1988).
\item[{[H2]}]   D. Happel, On Gorenstein algebras. In: Representation theory of finite groups and
finite-dimensional algebras (Proc. Conf. at Bielefeld, 1991). Progress in Math, vol. 95,
pp. 389-404. Birkh\"{a}user, Basel (1991).
\item[{[IYa]}]  O. Iyama and D. Yang, Silting reduction and Calabi¨CYau reduction, arXiv: 1408.2678 (2014).
\item[{[IYo]}]  O. Iyama and Y. Yoshino, Mutations in triangulated categories and rigid Cohen-Macaulay
modules, Invent. Math. 172 (2008), 117-168.
\item[{[KV]}]  B. Keller and D. Vossieck, Aisles in derived categories, Bull. Soc. Math. Belg. S¨¦r. A 40 (2) (1988),
239-253.
\item[{[KZ]}]  S. K\"{o}nig and A. Zimmermann, Derived equivalences for group rings, Lect. Note Math. 1685 (1998).
\item[{[K]}]  H.  Krause, The stable derived category of a noetherian scheme, Compos. Math.
141 (2005), 1128-1162.
%%
%\item[{[K2]}]  H. Krause, Cohomological length functions, arXiv: 1209.0540 (2012).
%
\item[{[MT]}]  N. Mahdou and M. Tamekkante, Mediterr. J. Math. 8 (2011), 293-305.
\item[{[O]}] D.  Orlov,  Triangulated categories of singularities and D-branes in Landau-Ginzburg
models, Proc. Steklov Inst. Math. 246 (3) (2004), 227-248.
\item[{[RS]}]   J. Rada and M. Saorin, Rings characterized by (pre)envelopes and (pre)covers of
their modules, Comm. Algebra 26 (1998), 899-912. %
\item[{[R1]}]   J. Rickard, Morita theory for derived categories, J. Lond. Math. Soc. 39 (2) (1989), 436-456. %
\item[{[R2]}] J. Rickard, Derived categories and stable equivalence. J. Pure Appl. Algebra 61 (1989), 303-317.
\item[{[W1]}] J. Wei, $\omega$-Gorenstein modules, Comm. Algebra 36 (5) (2008), 1817-1829.
\item[{[W2]}] J. Wei, Semi-tilting complexes, Israel J. Math. 194 (2013), 871-893.
\end{itemize}
%\end{thebibliography}
}

\end{document}